\documentclass[graybox]{svmult}

\usepackage{mathptmx}       
\usepackage{helvet}         
\usepackage{courier}        
\usepackage{type1cm}        
\usepackage{makeidx}         
\usepackage{graphicx}        
\usepackage{multicol}        
\usepackage[bottom]{footmisc}
\usepackage{amssymb}

\usepackage{hyperref}
\hypersetup{
    pdftitle={On Some Distributed Disorder Detection}, 
    pdfauthor=Krzysztof Szajowski, 
    bookmarks=true,
    colorlinks,
    urlcolor=blue,
    filecolor=magenta,
    citecolor=green, 
    linkbordercolor={1 1 1}, 
    citebordercolor={1 1 1},  
    urlbordercolor={ 1 1 1}  
  } 

\def\bbE{{\mathbb E}}
    \def\BbbE{\mathbb E}
    \def\BbbN{\mathbb N}

\def\bP{{\bf P}}
\def\bE{{\bf E}}

\def\cC{{\mathcal C}} 
\def\cF{{\mathcal F}}
\def\cW{{\mathcal W}}
\def\cL{{\mathcal L}}

\def\frt{{\mathfrak t}}
\def\frF{{\mathfrak F}}
\def\frS{{\mathfrak S}}
\def\frN{{\mathfrak N}}

\def\one{{\mathbb I}}

       \def\vecx{\overrightarrow{x}}
       \def\vecX{\overrightarrow{X}}

\toctitle{Distributed Disorder Detection}
\tocauthor{K.Szajowski}
\begin{document}

\title*{On Some Distributed Disorder Detection}

\titlerunning{Distributed Disorder Detection}

\author{Krzysztof Szajowski%
\thanks{The author would like to thank the Isaac Newton Institute for Mathematical Sciences, Cambridge, for support
and hospitality during the programme "Inference for Change-Point and Related Processes", where a part of the work on this paper was undertaken.} 
}
\authorrunning{K.Szajowski}
\institute{Krzysztof Szajowski \at Inst. of Math. and CS, Wroc\l{}aw Univ. of Tech., Wyb. Wyspia\'{n}skiego 27, PL-50-370 Wroc\l{}aw; \email{Krzysztof.Szajowski@pwr.edu.pl}
\url{http://im.pwr.edu.pl/\~szajow}}
%
%
\maketitle
\vspace{-2.79cm}
\abstract{Multivariate data sources with components of different information value seem to appear frequently in practice. Models in which the components change their homogeneity at different times are of significant importance. The fact whether any changes are influential for the whole process is determined not only by the moments of the change, but  also depends on which coordinates. This is particularly important in issues such as reliability analysis of complex systems and the location of an intruder in surveillance systems. In this paper we developed a mathematical model for such sources of signals with discrete time having the Markov property given the times of change. The research also comprises a multivariate detection of the transition probabilities changes at certain sensitivity level in the multidimensional process. Additionally, the observation of the random vector is depicted. Each chosen coordinate forms the Markov process with different transition probabilities before and after some unknown moment. The aim of statisticians is to estimate the moments based on the observation of the process. The Bayesian approach is used with the risk function depending on measure of chance of a false alarm and some cost of overestimation.  The moment of the system's disorder is determined by the detection of transition probabilities changes  at some coordinates. The overall modeling of the critical coordinates is based on the simple game. 
}
\vspace{-3ex}
\section{\label{a50intro}Introduction}
\vspace{-1.2ex}
The aim of the study is to investigate the mathematical model of a multivariate surveillance system introduced in \cite{sza11:multivariate}. in the model there is net $\frN$ of $p$ nodes. At each node the state is the signal at moment $n\in \BbbN$ which is at least one coordinate of the vector $\vecx_n\in\bbE\subset\Re^m$. The distribution of the signal at each node has two forms that depend on the state of surrounding. The state of the system changes dynamically. We consider the discrete time  signal observed as $m\geq p$ dimensional process on the probability space $(\Omega,\cF,\bP)$. The Markov processes, which are observed at each node,  are non homogeneous with two homogeneous segments as they have different transition probabilities (see \cite{sarsza11:transition} for details). The visual consequence of the transition distribution changes at moment $\theta_i$, $i\in\frN$ is a change of its character. In order to avoid false alarm the confirmation from other nodes is needed. The family of subsets (coalitions) of nodes is defined in such a way that the decision of all members of a given coalition is equivalent to the claim  that the disorder appeared in the net. It is not certain, however that the disorder has taken place. The aim is to define the rules of nodes and a construction of the net decision based on individual nodes claims. Various approaches can be found in the recent research that refer to the  description of such systems (see e.g. \cite{tarvee08:sensor}, \cite{ragvee10:Markov}). The problem is quite similar to a pattern recognition with multiple algorithm when the results of fusions of individual algorithms are unified to a final decision. In the study two different approaches are proposed. Both are based on the simple game defined on the nodes. The naive methods determine the system disordering by fusion individual node strategies. This construction of the individual decisions is based on the observation at each node separately.

The advanced solution of Bayesian version of the multivariate detection with a common fusion center is based on a stopping game defined by a simple game related to the observed signals. The individual decisions are based the analysis of the processes observed at all nodes and knowledge of nodes' interaction (the simple game). The  sensors' strategies are constructed as an equilibrium strategy in a non-cooperative stopping game with a logical function defined by a simple game (which aggregates their decision). 

The general description of such multivariate stopping games has been formulated by Kurano, Yasuda and Nakagami in the case when the aggregation function is defined by the voting majority rule \cite{kyn} and the observed sequences of the random variables are independent, identically distributed. It was Ferguson~\cite{fer05:MR2104375} who substituted the voting aggregation rules by a simple game. The Markov sequences have been investigated by the author and Yasuda~\cite{szayas95:voting}.

The model of detection of the disorders at each sensor is presented in the next section. It allows to define the individual payoffs of the players (sensors). 
The final decision based on the state of the sensors is given by the fusion center and it is described in Section~\ref{a50strategiesOFsensors}. The natural direction of further research is formulated in the same section. 
\vspace{-4ex}

\section{\label{a50disorderONsensor} Detection of disorder at sensors}
\vspace{-1.2ex}
Following the consideration presented in Section~\ref{a50intro}, let us suppose that the process $\{\vecX_n,n\in\BbbN\}$, $\BbbN=\{0,1,2,\ldots\}$, is observed sequentially in such a way that each sensor, \emph{e.g.} $r$th one gets its coordinates in the vector $\vecX_n$ at moment $n$. By assumption, it is a stochastic sequence that has the Markovian structure which is given random moment $\theta_r$ in such a way that the process after $\theta_r$ starts from state $\vecX_{\theta_r-1} $. The objective is to detect these moments based on the observation of $\vecX_{n}$ at each sensor separately. There are some results on the discrete time case of such disorder detection which generalize the basic problem stated by Shiryaev~in~\cite{shi61:detection} (see e.g. Brodsky and Darkhovsky~\cite{brodar93:nonparametr}, Bojdecki~\cite{boj79:disorder}) in various directions. In the early papers the observed sequence has independent elements given disorder moment. The sequences with dependent observations are subject of investigation by  Yoshida~\cite{yos83:complicated}, Szajowski~\cite{sza92:detection}, Yakir~\cite{yak94:finite}, Moustakides~\cite{Mou98:DisProcess} and Mei~\cite{Mei06:Yakir}.   

The application of the model for the detection of traffic anomalies in networks was discussed by Tartakovsky et al.~\cite{tarroz06:intrusions}. The version of the problem when the moment of disorder is detected with given precision will be used here (see~\cite{sarsza11:transition}).
\vspace{-.8cm}

\subsection{\label{a50a50sformProblem} Formulation of the problem}
\vspace{-.5cm}
The observable random variables $\{\vecX_n\}_{n \in \BbbN}$ are consistent with the filtration $\mathcal{F}_n$ (or $\cF_n = \sigma(\vecX_0,\vecX_1,\ldots,\vecX_n)$). The random vectors $\vecX_n$ take values in $(\bbE, \mathcal{B})$, where $\bbE\subset\Re^m$. On the same probability space there are defined unobservable (hence not measurable with respect to $\cF_n$) random variables $\{\theta_r\}_{r=1}^m$ which have the following geometric distributions:
\begin{eqnarray*}
\label{a50a50rozkladyTeta}
\bP(\theta_r = j) &=&\pi_r\one_{\{j=0\}}(j)+(1-\one_{\{j=0\}}(j))(1-\pi_r) p_r^{j-1}q_r, 
\end{eqnarray*}
where $\pi_r,q_r=1-p_r \in (0,1)$, $j=0,1,2,\ldots$.

The sensor $r$ follows the process which is based on switching between two, time homogeneous and independent the Markov processes $\{X_{rn}^i\}_{n \in \BbbN}$, $i=0,1$, $r\in\frN$ with the state space $(\bbE, \mathcal{B})$. These are both independent of $\{\theta_r\}_{r=1}^m$. Moreover, the processes $\{X_{rn}^i\}_{n \in \BbbN}$ have transition densities 
\begin{eqnarray*}\label{a50a50TransProbab}
\bP_x^{i}(X_{r1}^{i}\in B)&=&\bP(X_{r1}^{i}\in B|X_{r0}^{i}=x)=\int_Bf_x^{ri}(y)\mu(dy).
\end{eqnarray*}
The random processes $\{X_{rn}\}$, $\{X_{rn}^0\}$, $\{X_{rn}^1\}$ and the random variables $\theta_r$ are connected via the rule: $X_{rn}=X_{rn}^0\one_{\{n:n<k\}}(n)+ X_{r\;n+1-k}^1\one_{\{n:n\leq k\}}(n)$ on $\theta_r = k$, 
where $\{X_{rn}^1\}$ starts from $X_{r\;k-1}^0$ (but is otherwise independent of $X_{r\;\cdot}^0$).

For any $x\in \bbE$, $\pi_r\in [0,1]$, $c\in\Re_{+}$ and $\tau_r\in \frS^X$, where $\frS^X$ denotes the set of all stopping times with respect to the filtration
$\{\mathcal{F}_n\}_{n \in \BbbN}$, the risk associated with $\tau_r$ is defined as follows
$ \rho_r(x,\pi_r,\tau_r)= \bP_{x\;\pi_r}(   \tau_r  < \theta_r  ) + c\bE_{x\;\pi_r}\max\{\tau_r-\theta_r,0\},
$
where $\bP_{x\;\pi_r}(\tau_r<\theta_r)$ is the probability of false alarm and  $\bE_{x\;\pi_r}\max\{\tau_r-\theta_r,0\}$ is the average delay of detecting correctly the occurrence of disruption.

Every sensor is looking for the stopping time $\tau_r^{*}\in \frS^X$ such that for every $(x\,;\pi_r)\in\bbE\times [0,1]$
\begin{equation}
\label{a50PojRozregCiagowMark-Problem}
\rho^\star(x,\pi_r)=\rho_r(x,\pi_r,\tau_r^\star)=\inf_{\tau_r\in \frS^X}\rho_r(x,\pi_r,\tau_r).
\end{equation}
\vspace{-.75cm}

\subsection{The optimal detection problem as an optimal stopping problem}
\vspace{-.5cm}
In case of the independent sequence given the disorder moment  the construction of $\tau^{*}$ through the transformation of the problem to the optimal stopping problem for the Markov process $(X_n,\Pi^{\pi_r}_{r\; n})$ can be made where $\Pi^{\pi_r}_{r\; n}$ is the posterior process (see e.g. \cite{shi78:book}). It is stated that $\Pi^{\pi_r}_{r0} = \pi_r$,  $\Pi^{\pi_r}_{rn} = \bP^{\pi_r}\left(\theta_r \leq n \mid \mathcal{F}_n\right)$, for $n = 1, 2, \ldots$, is designed as information about the distribution of the disorder instant $\theta_r$. Moreover, 
\begin{equation}\label{a50payoff2}
\rho_r(x,\pi_r,\tau_r)= \bE^{x,\pi_r}\left\{(1-\Pi^{\pi_r}_{r\tau_r})+c\sum_{k=0}^{\tau_r-1}\Pi^{\pi_r}_{r\;k}\right\}.
\end{equation}
The family of the Markov random functions $\{\Pi^{\pi_r}, \pi_r\in [0,1]\}$ can be associated with a Markov process with discrete time $\Pi=(\pi_n,\cF_n,\bP^{\pi_r})$, for $n\geq 0$, having the same transition probabilities as each Markov random function  is presented as $\Pi^{\pi_r}$, $\pi_r\in[0,1]$. 
\vspace{-.75cm}

\subsection{\label{a50OSProblem1} The optimal stopping problem with observation costs}
\vspace{-.5cm}
The problem of minimization of the risk (\ref{a50payoff2}) can be solved as the special optimal stopping problem.
As it is shown in \cite{pesshi06:optimal}, p.22-23, the problem can be transformed to the optimal stopping problem for the time-homogeneous two dimensional Markov chain without observation costs.  The Wald-Bellman equation which solves (\ref{a50PojRozregCiagowMark-Problem}) takes the form:
\begin{equation}\label{a50WBeq2}
\rho^\star(x,\pi_r)=\min\{1-\pi_r,c\pi_r+\bE_{x,\pi_r}\rho^\star(x_1,\pi_1)\}.
\end{equation}

\vspace{-3ex}
\section{\label{a50coopnobcoop} The aggregated decision via the cooperative game}
\vspace{-1.2ex}
There are various methods combining the decisions of several classifiers or sensors. The methods based on winning coalitions in the simple game presented in \cite{sza11:multivariate} will be used. The obvious changes are the consequence that in the model considered now the aim is to minimize the risk. We apply two methods of decision aggregation. In the first one, based on the optimal disorder detection strategies, we apply the aggregation method. This approach does not guarantee that the obtained system disorder detection will have certain stability or equilibrium properties. 

In the second approach each ensemble member contributes to some degree to the decision at any point of the sequentially delivered states. The fusion algorithm takes into account all the decision outputs from each ensemble member and comes up with an ensemble decision in such a way that the solution is an equilibrium point in an  antagonistic, no-zero sum game.  
\vspace{-.6cm}

\subsection{A simple game}
\vspace{-.5cm}
Let us assume that there are many nodes which absorb information and make decisions if the disorder has appeared or not. The final decision is made in the fusion center which aggregates the information from all sensors.  

The voting decision is made according to the rules of \emph{a simple game}. Let us recall that a coalition is a subset of the players.  Let ${\cC}=\{C:C\subset \frN\}$ denote the class of all coalitions.  
\emph{A simple game} (see \cite{Owe95:MR1355082}, \cite{fer05:MR2104375}) is a coalition game having the characteristic function of $\phi(\cdot):\cC\rightarrow\{0,1\}$. 

Let us denote $\cW=\{C\subset\frN:\phi(C)=1\}$ and ${\cL}=\{C\subset\frN:\phi(C)=0\}$. The coalitions in $\cW$ are called the winning coalitions, and those from $\cL$ are called the losing coalitions.
By assumption, the characteristic function satisfies the properties: $\frN\in\cW$; $\emptyset\in \cL$; (the monotonicity): $T\subset S\in \cL$ implies $T\in \cL$.
\vspace{-.5cm}

\subsection{The aggregated decision rule}
\vspace{-.5cm}
When the simple game is defined and the players can vote presence or absence, $x_i=1$ or $x_i=0$, $i\in\frN$ of the local disorder then the aggregated decision is given by the logical function
\begin{equation}\label{a50aggregateFUNCTION}
\delta(x_1,x_2,\ldots,x_p)=\sum_{C\in\cW}\prod_{i\in C}x_i\prod_{i\notin C}(1-x_i).
\end{equation}
For the logical function $\delta $ we have (cf \cite{kyn}) 
\[
\delta (x^1,\ldots ,x^p)=
x^i\cdot \delta (x^1,\ldots ,\stackrel{i}{\breve{1}}%
,\ldots,x^p)
+\overline{x}^i\cdot \delta (x^1,\ldots ,\stackrel{i}{\breve{0}}%
,\ldots ,x^p).
\]
\vspace{-1.5cm}

\subsection{\label{a50agregateSENSORS} Aggregated sensors strategies}
\vspace{-.5cm}
For any stopping times $\{\tau_i\}_{i=1}^p$ with respect of the filtration $\{\mathcal{F}_n\}_{n \in \BbbN}$ we have the representation by the individual stopping strategies $\sigma_n^i(\tau)=\one_{\{\omega:\tau_i\geq n\}}$. The aggregate function applied to the individual stopping times will construct the detection strategy $\sigma_n$ of  the system disorder. The stopping time from the individual stopping strategy is constructed as $\tau=\inf\{0\leq n\leq N:\sigma_n\prod_{k=1}^{n-1}(1-\sigma_k)=1\}$. 

This aggregation method is the basement of both constructions. In the \emph{naive} algorithm it is applied to the optimal individual strategies of the sensors constructed as the solution of the optimal stopping problem (\ref{a50PojRozregCiagowMark-Problem}). 

In the multivariate stopping game approach the aggregation of the individual decision is used to construct the set of admissible strategies. The details are the subject of the next section.  

\vspace{-1.3ex}
\vspace{-.5cm}
\section{\label{a50noncooperative}A non-cooperative detection problem}
\vspace{-1.2ex}
Following the results of the author and Yasuda~\cite{szayas95:voting} the multilateral stopping of a Markov chain problem can be described in the terms of the notation used in the non-cooperative game theory (see \cite{nas51:noncoop}, \cite{Owe95:MR1355082}). This approach can be applied to the distributed disorder detection by reformulation of the problem to the multilateral stopping problem. The important issue is the representation of the expected risk in the disorder detection problem for one sensor given in (\ref{a50payoff2}). 

Let us denote $\sigma ^i=(\sigma _1^i,\ldots ,\sigma _N^i)$ and let $\frS^i$ be the set of ISSs of player $i$, $i=1,2,\ldots ,p$ (see \cite{kyn}). Define $\frS=\frS^1\times \ldots\times \frS^p$ the set of the stopping strategy (SS). The factual stopping of the observation process (the estimate of the system disorder moment), and the players realization of the payoffs are defined by the stopping strategy exploiting $p$-variate logical function $\delta :\{0,1\}^p\rightarrow \{0,1\}$. Since $\delta $ is fixed during the analysis we write $\frt(\sigma )=\frt_\delta (\sigma )$. 

We have $\{\omega \in \Omega : \frt_\delta (\sigma )=n\} =\bigcap\nolimits_{k=1}^{n-1}\{\omega \in \Omega : \delta (\sigma _k^1,\sigma_k^2,\ldots,\sigma _k^p)=0\}\cap \{\omega \in \Omega :\delta (\sigma_n^1,\sigma _n^2,\ldots,\sigma _n^p)=1\}\in \frF_n$, then the random  variable $\frt_\delta (\sigma )$ is the stopping time with respect to $\{\frF_n\}_{n=1}^N$. For any stopping time $\frt_\delta (\sigma )$ and $i\in \{1,2,\ldots ,p\}$, let $\rho_i(X_{\frt_\delta (\sigma )},\Pi_{\frt_\delta (\sigma )},\delta (\sigma ))=\rho_i(X_n,\Pi_n,n)\one_{\{\frt_\delta (\sigma )=n\}}+\limsup_{n\rightarrow \infty }\rho_i(X_n,\Pi_n,n)\one_{\{\frt_\delta (\sigma )=\infty \}}$.  
(cf \cite{shi78:book}, \cite{szayas95:voting}). If players use SS $\sigma \in \frS$ and the individual preferences are converted to the effective stopping time by the aggregate rule $\delta $, then player $i$ gets $\rho_i(X_{\frt_\delta (\sigma )},\Pi_{\frt_\delta (\sigma )})$. 

Let ${}^{*}\!\sigma =({}^{*}\!\sigma ^1,\ldots ,{}^{*}\!\sigma ^p)^T\in\frS $ and ${}^{*}\!\sigma (i)=({}^{*}\!\sigma ^1,\ldots ,{}^{*}\!\sigma ^{i-1},\sigma ^i,{}^{*}\!\sigma ^{i+1},\ldots,{}^{*}\!\sigma^p)^T$.

\begin{definition}
\label{a50equdef}{\rm (cf. \cite{szayas95:voting})} For the fixed aggregate rule $\delta $ the strategy
${}^{*}\!\sigma \in \frS$ is an equilibrium strategy if for each $i\in \{1,2,\ldots ,p\}$ and any $\sigma^i\in \frS^i$ we have
\begin{equation}\label{a50defequ} 
\rho_i(x,\pi_i,\frt_\delta ({}^{*}\!\sigma))\leq \rho_i(x,\pi_i,\frt_\delta({}^{*}\!\sigma(i))). 
\end{equation}
\end{definition}

The set $\frS$, the vector of the utility functions $f=(f_1,f_2,\ldots, f_p)$ and the monotone rule $\delta $ define the non-cooperative game $\cal{G}$ = ($\frS$,$f$,$\delta$). The construction of the equilibrium strategy $ {}^{*}\!\sigma \in \frS$ in $\cal{G}$ is provided in \cite{szayas95:voting}. In the case of the considered distributed disorder detection problem we have $f_i(x,\pi)=1-\pi$.    

With each ISS of player $i$ the sequence of stopping events $D_n^i=\{\omega :\sigma _n^i=1\}$ is associated. For each aggregate rule $\delta$ there exists the corresponding set value function $\Delta :\frF\rightarrow \frF$ such that $\delta (\sigma_n^1,\ldots ,\sigma _n^p)= \delta \{\one_{D_n^1},\ldots,\one_{D_n^p}\}= \one_{\Delta(D_n^1,\ldots,D_n^p)}$. For the solution of the considered game the important class of ISS and the stopping events can be defined by subsets ${\it{C}}^i \in \cal{B}$ of the state space $\BbbE$. A given set ${\it{C}}^i\in\cal{B}$ will be called the stopping set for player $i$ at moment $n$ if $D_n^i= \{\omega :X_n\in {\it{C}}^i\}$ is the stopping event. 

Let $g_i$ be the real, integrable functions defined on $\BbbE\times [0,1]$ and ${\it{C}}^i\in \cal{B}$. Let ${}^i\!D_1(A)=\Delta (D_1^1,\ldots ,D_1^{i-1},A,D_1^{i+1},\ldots ,D_1^p)$. For fixed $D_n^j=\{\omega :X_n\in {\it{C}}^i\}$, $ j=1,\ldots ,p$, $j\neq i$ define $\psi ({\it{C}}^i)={\bf E}_x\left[(1-\Pi^i_1)\one_{{}^i\!D_1(D_1^i)}+ g_i(X_1, \Pi_1)\one_{\overline{{}^i\!D_1(D_1^i)}}\right]$.  
\begin{lemma} $\label{a50optimal}$
Let ${\it{C}}^j\in \cal{B}$, $j=1,2,\ldots,p$,  $j\neq i$, be fixed. Then the set ${}^{*}\!{\it{C}}^i=\{x\in \BbbE:\rho_i(x)-g_i(x)\leq 0\}\in \cal{B}$ is such that $\psi ({}^{*}\!{\it{C}}^i)=\inf\limits_{{\it{C}}^i\it{\in }\cal{B} }\psi ({\it{C}}^i) $
and
\begin{eqnarray}\label{a50optset}
\psi ({}^{*}\!{\it{C}}^i)
& = & {\bf E}_{x,\pi}(1-\Pi^i_1-g_i(X_1,\Pi_1))^{+}\one_{{}^{i}\!D_1(\emptyset )}\\
\nonumber& & - {\bf E}_{x,\pi}(1-\Pi^i_1-g_i(X_1,\Pi_1))^{-}\one_{{}^{i}\!D_1(\Omega )} +{\bf
E}_{x,\pi}g_i(X_1,\Pi_1). 
\end{eqnarray}
\end{lemma}

Based on Lemma \ref{a50optimal} we derive the recursive formulae defining the
equilibrium point and the equilibrium payoff for the finite horizon detection problem. 
\vspace{-.5cm}

\subsection{The finite horizon detection problem\label{a50finite}} 
\vspace{-.5cm}
Let horizon $N$ be finite and the equilibrium strategy ${}^{*}\!\sigma $ exist. We denote $\rho_{i,N}(x,\pi)={\bf E}_{x,\pi}\rho_i(X_{t({}^{*}\!\sigma )},\Pi_{t({}^{*}\!\sigma )})$ the equilibrium payoff of $i$-th player when $X_0=x$. Let $\frS_n^i=\{\{\sigma_k^i\},k=n,\ldots ,N\}$  and $\frS_n=\frS_n^1\times \frS_n^2\times \ldots \times\frS_n^p$. 

Denote $t_n=t_n(\sigma )=t(^n\sigma )=\inf \{n\leq k\leq N:\delta (\sigma _k^1,\sigma_k^2,\ldots ,\sigma _k^p)=1\}$ to be the stopping time not earlier than $n$. 
Let $\rho_{i,N-n+1}(X_{n-1},\Pi_{n-1})=\rho_i(X_{t_n({}^{*}\!\sigma )},\Pi_{t_n({}^{*}\!\sigma )},t_n({}^{*}\!\sigma ))$. At $n=N$ we have $\rho_{i,0}(x,\pi)=\rho_i(x,\pi,N)$. Let us assume that the process is not stopped up to moment $n$ and the players are using the equilibrium strategies ${}^{*}\!\sigma _k^i$, $i=1,2,\ldots ,p,$ at $k=n+1,\ldots ,N$. Choose player $i$ and assume that other players are using the equilibrium strategies ${}^{*}\!\sigma _n^j$, $j\neq i$, and player $i$ is using strategy $\sigma_n^i$ defined by the stopping set ${\it{C}}^i$. Then the expected payoff $\varphi_{N-n}(X_{n-1},{\it{C}}^i) $ of player $i$ in the game starting at $n$, when the state of a Markov chain at $n-1$ is $X_{n-1\mbox{,}}$ is equal to
\[
 \varphi _{N-n}(X_{n-1},\Pi_{n-1},{\it{C}}^i)=
 {\bf E}_{X_{n-1},\Pi_{n-1}}\left[(1-\Pi_n^i)\one_{{}^{i*}\!D_n(D_n^i)}+ 
 \rho_{i,N-n}(X_n,\Pi_n)
 \one_{\overline{{}^{i*}\!D_n(D_n^i)}}\right], 
\]
where
${}^{i*}\!D_n(A)=\Delta({}^{*}\!D_n^1,\ldots,{}^{*}\!D_n^{i-1},A, 
{}^{*}\!D_n^{i+1},\ldots,{}^{*}\!D_n^p)$. 

By Lemma \ref{a50optimal} the conditional expected gain $\varphi _{N-n}(X_{N-n}, {\it{C}}^i)$ attains the maximum on the stopping set ${}^{*}\!{\it{C}}_n^i=\{x\in \bbE:f_i(x)-v_{i,N-n}(x)\leq 0\}$ and 
\setcounter{equation}{0}
\small
\begin{eqnarray*} \label{a50valueatn}
v_{i,N-n+1}&(&X_{n-1},\Pi_{n-1})-c_i \Pi_{n-1}^i ={\bf E}_x[(1-\Pi^i_n-v_{i,N-n}(X_n,\Pi_n))^{+} 
\one_{{}^{i*}\!D_n(\emptyset)}|\frF_{n-1}] \\
\nonumber& &- {\bf E}_x[(1-\Pi^i_n-v_{i,N-n}(X_n))^{-} \one_{{}^{i*}\!D_n(\Omega)}|\frF_{n-1}] + {\bf E}_x[v_{i,N-n}(X_n,\Pi_n)|\frF_{n-1}] 
\end{eqnarray*}\normalsize
${\bf P}_x-$a.e..
This reasoning allows to formulate the following construction of the equilibrium strategy and the equilibrium value for the game $\cal{G}$.

\begin{theorem}
In the game $\cal{G}$ with finite horizon $N$ we have the following solution. 
\begin{description}
\item[(i)] The equilibrium value $v_i(x,\pi)$, $i=1,2,\ldots ,p$, of the game 
$\cal{G}$ can be calculated recursively as follows: 
$v_{i,0}(x,\pi)=1-\pi_i$ and
for $n=1,2,\ldots ,N$, $i=1,2,\ldots ,p$ we have ${\bf P}_x-$a.e. 
\small
\begin{eqnarray*}
v_{i,n}(X_{N-n},\Pi_{N-n})&-&c_i \Pi^i_{N-n}={\bf E}_{x,\pi}[v_{i,n-1}(X_{N-n+1},\Pi_{N-n+1})|\frF_{N-n}]\\
&+&{\bf E}_{x,\pi}[(1-\Pi^i_{N-n+1})-v_{i,n-1}(X_{N-n+1}),\Pi_{N-n+1})^{+} \one_{{}^{i*}\!D_{N-n+1}(\Omega )}|\frF_{N-n}] \\ 
&-& {\bf E}_{x,\pi}[((1-\Pi^i_{N-n+1})-v_{i,n-1}(X_{N-n+1},\Pi_{N-n+1}))^{-} 
\one_{{}^{i*}\!D_{N-n+1}(\emptyset )}|\frF_{N-n}]. 
\end{eqnarray*}
\normalsize
\item[(ii)] The equilibrium strategy ${}^{*}\!\sigma \in \frS$ is defined by the SS of the players ${}^{*}\!\sigma _n^i$, where ${}^{*}\!\sigma _n^i=1$ if $ X_n\in {}^{*}\!{\it{C}}_n^i$, and ${}^{*}\!{\it{C}}_n^i=\{x\in \bbE : f_i(x)-v_{i,N-n}(x) \leq 0\}$, $n=0,1,\ldots ,N$.
\end{description}

We have $v_i(x,\pi)=v_{i,N}(x,\pi)$, and ${\bf E}_{x,\pi}(1-\Pi^i_{t({}^{*}\!\sigma )})=v_{i,N}(x)$, $i=1,2,\ldots ,p$.
\end{theorem}
\nopagebreak



\vspace{-3ex}
\section{\label{a50strategiesOFsensors}Determining the strategies of sensors}
\vspace{-1.1ex}
Based on the model constructed in Sections~\ref{a50disorderONsensor}--\ref{a50noncooperative} for the net of sensors with the fusion center determined by a simple game, one can determine the rational decisions of each nodes. The rationality of such a construction refers to the individual aspiration for the highest sensitivity to detect the disorder without a false alarm. The Nash equilibrium fulfills the requirement that nobody deviates from the equilibrium strategy, otherwise its expected risk will be higher. 

The proposed model disregards the correlation of the signals. It is also assumed that the fusion center has complete information about the signals and that the information is available at each node. The method of a cooperative game was used in \cite{ghakri10:cooperativeMR2780188} to find the best coalition of sensors in the problem of the target localization. The approach which is proposed in the study shows the possibility of modeling the detection problem by multiple agents at a general level. 
\vspace{-2.3ex}



\end{document}